\documentclass{amsart}
\usepackage{amssymb,latexsym,amsmath}

\theoremstyle{plain}
\newtheorem{theorem}{Theorem}[section]
\newtheorem{corollary}[theorem]{Corollary}
\newtheorem{lemma}[theorem]{Lemma}
\newtheorem{proposition}[theorem]{Proposition}
\newtheorem{definition}[theorem]{Definition}

\newtheorem*{remark}{Remark}

\newcommand{\half}{\frac{1}{2}}

\numberwithin{equation}{section}

\font\cmcsc=cmcsc10 at 8pt

\begin{document}

\title[geometric brownian motion]
        {exponential martingales and time integrals of brownian motion}
\author{Victor Goodman}
\address{Mathematics Department\\
         Indiana University\\
         Bloomington, IN 47405}
\email{goodmanv@indiana.edu}
\author{Kyounghee Kim}
\address{Mathematics Department\\
         Florida State University\\
         Tallahassee, FL  32306-4510}
\email{kim@math.fsu.edu}

\maketitle

\begin{abstract}
We find a simple expression for the probability density of  \\$\int \exp (B_s - s/2) ds$ in terms of its distribution function and the distribution function for the time integral of $\exp (B_s + s/2)$.  The relation is obtained with a change of measure argument where expectations over events determined by the time integral are replaced by expectations over the entire probability space.  We develop precise information concerning the lower tail probabilities for these random variables as well as for time integrals of geometric Brownian motion with arbitrary constant drift.  In particular,  
$E[\ \exp\big(\theta / \int \exp (B_s )ds\big) ]$ is finite iff $\theta < 2$.  We present a new formula for the price of an Asian call option.\hfill\break
\phantom{as}\hfill\break
\noindent {\cmcsc Keywords}:\ Girsanov theorem, Geometric Brownian Motion, Asian option.\hfill\break
\noindent {\cmcsc Subject Classification}:\   Primary 60J65, 60H30\ \      Secondary 91B28
\end{abstract}

\vspace{4ex}

\section{Introduction}

Time integrals of one-dimensional geometric Brownian motion have appeared in financial models where certain expected values are the computed prices of  `Asian options'.  Approximate values for some of these expectations were obtained in \cite{RS}.  Some fundamental work in \cite{Y92} and \cite{GY} focuses on the distribution and density functions of time integrals.  Matsumoto and Yor [MY 1--3] studied time-integrals of the form
$$\int_0^t \exp(-2B_s - \mu s)ds$$
for $\mu > 0$ and obtained some interesting identities in law (Theorem 1.1 and 2.3 in \cite{MY3}).  In a recent paper Bertoin and Yor \cite{BY} considered the time-integral of Brownian functions over an infinite time interval.  The distributions of time integrals over an infinite interval are known \cite{D2}:
$$\int_0^\infty \exp(-2B_s - \mu s)ds = \frac{1}{2\gamma_{\mu/2}}$$
where the law of  $\gamma_{\mu/2}$ is a gamma distribution with index $\mu/2$.

We present new relationships between the density functions and distribution functions for these random variables.  Our results provide a precise description of the lower tail of these distributions and we settle several moment questions involving an exponent which is the reciprocal of a time integral.  We also consider other  moments where the exponent includes further exponential terms involving Brownian motion. Certain of these double exponential terms have a finite expected value.

We let $B_t$ denote a one-dimensional Brownian motion. Our notation for a time integral of exponential Brownian motion is somewhat nonstandard. We let $M_t$ denote the simple exponential martingale
$$M_t = \exp(B_t - \frac{t}{2})$$
and we define its time integral as

$$A_t = \int_0^t M_sds.$$

It is known (see Dufresne, Corollary 2.3, \cite{D}) that for $\alpha >0$ and $\beta > \half$

$$E[\exp(\frac{M_t}{\alpha + \beta A_t})] < \infty$$

\noindent  We show that this moment, in fact, is finite for $\beta = \half$.

In [Y92,02] and \cite{BTW} the authors use $A_t$ to denote the time integral of exponential Brownian motion without drift, or an integral with a drift other than  $ -1/2$.  We obtain some distributional properties of the random  variable $A_t$.  In particular we present a formula for its probability density and we obtain a sharp lower tail estimate for its distribution.

Our method of argument uses a delicate change of measure, where the process

\begin{equation}\label{other} 
\tilde B_t := B_t - 2\log(1- \frac{y}{2}A_t)
\end{equation}

\noindent becomes a standard Brownian motion.  Here, $y$ is a positive constant.  Of course, this path translation has a singularity at the random time $\tau_\infty$ defined by
\begin{equation}\label{last} 
A_{\tau_\infty} = \frac{2}{y}
\end{equation}
Therefore, we will consider sample paths only for $t$ strictly less than $ \tau_\infty$.  Matsumoto and Yor (Section 2 in \cite{MY3}) considered a change of measure where

$$\tilde B_t^{(\mu)} = B_t^{(\mu)} - 2\log(1+\gamma_\mu A_t)$$

\noindent and the {\em random coefficient} $\gamma_\mu $ is independent of the Brownian motion and has a gamma distribution with index $\mu$.  In contrast, our process $\tilde B_t$ is defined without this random factor and our process exists only up to the explosion time $\tau_\infty$.

\vspace{3ex}

\section{ The Change of Measure}

\vspace{3ex}

We follow the Girsanov formalism: A sufficient condition for a process
$$\tilde B_t := B_t - \int_0^t \theta_s ds$$
to be a Brownian motion over a compact time interval $[0,T]$ is that the process
\begin{equation}\label{girsan} 
\Lambda_t := \exp\big(  \int_0^{ t} \theta_s dB - \frac{1}{2}  \int_0^{  t}\theta_s^2ds
 \big)
\end{equation}
be a martingale over the time interval (see [KS]).   $ \tilde B_t $ is a standard Brownian motion w.r.t. the measure $Q$ defined by 
\begin{equation}\label{qqq} 
\frac{dQ}{dP} =\Lambda_T
\end{equation}
In view of equation (\ref{other}), we consider as our choice of $\theta_t$ the process
\begin{equation}\label{rrr} 
R_t =2 \frac{d}{dt}\log(1- \frac{y}{2}A_t) = -\frac{M_t}{y^{-1}  -   \frac{1}{2}A_t}
\end{equation}
Notice that $R_0 = -y$ and that $R_t$ is defined up to the random time $ \tau_\infty$.  The process $R_t$ has the following convenient and remarkable property:

\vspace{3ex}

%begin the lemma on the sde of R

\begin{lemma}\label{sde}  The process given in equation (\ref{rrr}) satisfies the SDE
\begin{equation}\label{drdr}
dR_t =  R_t dB_t  -  \frac{1}{2} R_t^2dt
\end{equation}
for $t < \tau_\infty$.
\end{lemma}

\vspace{2ex}

\begin{proof}  This is a simple calculation using the fact that $dM_t = M_tdB_t$.
\end{proof}

\vspace{3ex}

\begin{remark}
By writing equation (\ref{drdr}) in its integral form,
$$
R_t =-y+  \int_0^t R_s dB - \frac{1}{2}  \int_0^tR_s^2ds,
$$
we see that $R_t$  is essentially the exponent of the Girsanov density process it generates.  This unusual property of $R_t$ allows us to analyze the behavior of $A_t$ through a change of measure. 
\end{remark}
\vspace{3ex}

\begin{definition} For each $n = 1,2, \dots$ let $\tau_n$ denote the stopping time given by
$$\tau_n = \inf \{t: R_t \leq -n \}$$
\end{definition}

\vspace{2ex}

\noindent Although each stopping time, and $\tau_\infty$ as well,  depends on the choice of $y$ , we will omit mentioning their dependence on this parameter unless we explicitly change the value of $y$.  We see from equation (\ref{last}) and the fact that $ R_t$ is negative that 
$$\tau_n < \tau_\infty
\hbox{\hskip.2in and also \hskip.2in }
\lim_{n\to\infty} \tau_n = \tau_\infty  \hbox{\hskip.4in a. s.}$$

%proposition tilde begins here

\vspace{3ex}

\begin{proposition}\label{tilde} For each $n = 1,2, \dots$ the process
\begin{equation}\label{mart}
\Lambda_t^{(n)}  =   \exp[\ y+  R_{t\wedge \tau_n}  \ ]
\end{equation}
forms a martingale.  Moreover, the process

$$ B_t - 2\log(1- \frac{y}{2}A_{t\wedge \tau_n} )$$

\noindent is a standard Brownian motion for $t \leq T$ w.r.t. the probability measure

$$dQ = \Lambda_T^{(n)} dP$$

\end{proposition}

\vspace{3ex}

\begin{proof}  We choose a Girsanov density process as in equation (\ref{girsan}) by setting 

$$\theta_s = R_s 1_{\{ s< \tau_n\}}$$
\vskip.1in

With this choice, $\theta_s$ is a bounded adapted process and hence satisfies a Novikov condition (see Corollary 5.13 of \cite{KS}).  The Novikov condition is sufficent for $\Lambda_t$ to be a martingale.  In our case,
$$\Lambda_t =  \exp\big(  \int_0^{ t\wedge  \tau_n} R_s dB - \frac{1}{2}  \int_0^{  t\wedge \tau_n}R_s^2ds\big)$$
It follows from Lemma \ref{sde} that the exponent above is precisely $R_{ t\wedge  \tau_n} + y$.  That is,

$$\Lambda_t  = \exp(y+ R_{ t\wedge  \tau_n} )$$
\vskip.2in

This proves the martingale assertion of the proposition.  Moreover, the calculation in equation (\ref{rrr}) shows that 
$$\int_0^t \theta_s ds =   2\log(1- \frac{y}{2}A_{t\wedge \tau_n} )$$
and so the Girsanov theorem implies that the process
$$ B_t - 2\log(1- \frac{y}{2}A_{t\wedge \tau_n} )$$
is a standard Brownian motion on compact time intervals with the change of measure given by $\Lambda_T$.

\end{proof}

\vspace{3ex}

\section{The Correspondence Between $B_t$ and $\tilde B_t$}

\vspace{3ex}

\begin{definition}\label{dict}  For fixed $n= 1,2, \dots$ and $y > 0$ we let $\tilde B_t$ denote the process appearing in Proposition \ref{tilde}.  That is

\begin{equation}\label{brown}
\tilde B_t :=  B_t - 2\log(1- \frac{y}{2}A_{t\wedge \tau_n} )
\end{equation}
\noindent We also define 

$$\tilde M_t = \exp(\tilde B_t - t/2)$$

and

$$\tilde A_t = \int_0^t \tilde M_s ds$$

\end{definition}

\vspace{3ex}

\noindent We note that all quantities in this definition depend on our choice for $n$ and $y$.

\vspace{3ex}
%  proposition dictionary begins here
\begin{proposition}\label{dictionary}

\begin{equation}\label{dictm}
\tilde M_t = \frac{M_t}{(1-\frac{y}{2}A_{t\wedge \tau_n})^2}
\end{equation}

\noindent and

\begin{equation}\label{dicta}
 1+\frac{y}{2}\tilde A_{t\wedge \tau_n} = \frac{1}{1-\frac{y}{2} A_{t\wedge \tau_n}}
\end{equation}

\noindent and

\begin{equation}\label{dicty}
R_{t\wedge \tau_n}= -\frac{\tilde M_{t\wedge \tau_n}}{y^{-1}+\frac{1}{2}\tilde A_{t\wedge \tau_n} }
\end{equation}

\end{proposition}

\vspace{3ex}

\begin{proof}  From Definition \ref{dict} we have

$$\tilde M_t = \exp( B_t - t/2 - 2\log(1- \frac{y}{2}A_{t\wedge \tau_n} ))$$
$$=\frac{M_t}{(1-\frac{y}{2}A_{t\wedge \tau_n})^2}$$
Now if $t < \tau_n$ then

$$\frac{d}{dt}(1+\frac{y}{2}\tilde A_t) =\frac{y}{2}\tilde M_t = \frac{yM_t}{2(1-\frac{y}{2}A_t)^2}
=\frac{d}{dt}\big(\frac{1}{1-\frac{y}{2}A_t}\big)$$

\noindent Hence

$$1+\frac{y}{2} \tilde A_t  = \frac{1}{1-\frac{y}{2}A_t}$$
Finally, for $t$ in this same range,

$$R_t=-\frac{M_t}{y^{-1}  -   \frac{1}{2}A_t} $$

$$ = -\frac{\tilde M_t({1  -   \frac{y}{2}A_t})^2}{y^{-1}  -   \frac{1}{2}A_t}$$
$$ = -y\tilde M_t({1  -   \frac{y}{2}A_t})$$
$$ = -\frac{y\tilde M_t}{  1+   \frac{y}{2}\tilde A_t}$$

\noindent These equalities hold up to the time $\tau_n$, and these are the assertions (\ref{dicta}) and (\ref{dicty}) in the proposition.
\end{proof}

\vspace{3ex}

%  proposition integral begins here
\begin{proposition}\label{integral}  If $f(x,z)$ is a  nonnegative Borel-measurable function  and $y>0$ then

\begin{equation}\label{expect}
E[ f(M_t, R_t)\ ; \  A_t < 2/y] 
\end{equation}

$$=E[ f(\frac{ M_t}{(1+\frac{y}{2} A_t)^2},  \frac{- M_t}{  y^{-1}+   \frac{1}{2} A_t})
\exp(\frac{ M_t}{  y^{-1}+   \frac{1}{2} A_t}-y)]$$

\end{proposition}

\vspace{3ex}

\begin{proof}  For fixed $n$ and $y>0$ we consider

$$E[ f(M_t, R_t)\ ; \  \tau_n >t ] $$
Now 

$$ f(M_t, R_t)1_{\{  \tau_n >t\}} =  f(M_t, R_t)\exp (-R_t - y)\exp (R_t +y)1_{\{  \tau_n >t\}}$$
\vskip.1in

\noindent Since this function vanishes for $t\geq \tau_n$, each time parameter may be replaced by  $t\wedge \tau_n $.  This allows us to apply Proposition \ref{tilde} where we take

$$\Lambda_t^{(n)}1_{\{  \tau_n >t\}} = \exp (R_t +y)1_{\{  \tau_n >t\}}$$ 

We obtain the identity

\begin{equation}\label{equal}
E[ f(M_t, R_t)\ ; \  \tau_n >t ] = E_Q [  f(M_t, R_t)\exp (-R_t - y)1_{\{  \tau_n >t\}} ]
\end{equation}

\vskip.3in

\noindent Each term in the r.h. expected value can be expressed in terms of the Brownian motion $\tilde B_t$. We use the identities in Proposition \ref{dictionary} to see that

$$f(M_t, R_t)\exp (-R_t - y)1_{\{  \tau_n >t\}}= $$
$$f(\tilde M_t (1-\frac{y}{2}A_{t})^2, -\frac{\tilde M_t}{y^{-1}+\frac{1}{2}\tilde A_t })
\exp (\frac{\tilde M_t}{y^{-1}+\frac{1}{2}\tilde A_t }- y)1_{\{  \tau_n >t\}}$$
$$=f(\frac{\tilde M_t}{(1+\frac{y}{2}\tilde A_t)^2 }, -\frac{\tilde M_t}{y^{-1}+\frac{1}{2}\tilde A_t })
\exp (\frac{\tilde M_t}{y^{-1}+\frac{1}{2}\tilde A_t }- y)1_{\{  \tau_n >t\}}$$
\noindent Moreover, the event $ \tau_n >t$ equals the event

 $$  \min_{s\leq  t} R_s > -n$$
which in turn equals the event
 $$  \max_{s\leq  t} \frac{\tilde M_s}{y^{-1}+\frac{1}{2}\tilde A_s } <  n$$

 \noindent The r.h. expected value in equation (\ref{equal}) is then
 
 $$E_Q [ f(\frac{\tilde M_t}{(1+\frac{y}{2}\tilde A_t)^2 }, -\frac{\tilde M_t}{y^{-1}+\frac{1}{2}\tilde A_t })
\exp (\frac{\tilde M_t}{y^{-1}+\frac{1}{2}\tilde A_t }- y)
1_{\{   \max_{s\leq  t} \frac{\tilde M_s}{y^{-1}+\frac{1}{2}\tilde A_s} <  n\}}]$$

But since the integrand is nonnegative we may take the limit as $n\to\infty$ and obtain the limiting value

$$E_Q [ f(\frac{\tilde M_t}{(1+\frac{y}{2}\tilde A_t)^2 }, -\frac{\tilde M_t}{y^{-1}+\frac{1}{2}\tilde A_t })
\exp (\frac{\tilde M_t}{y^{-1}+\frac{1}{2}\tilde A_t }- y)]$$ 
 
 \noindent In addition, since 
 
 $$ \lim_{n\to\infty} \tau_n = \tau_\infty$$
 as $ \tau_\infty$ is defined in equation (\ref{last})  we see that the limit of the l.h. side of equation (\ref{equal}) is 
 
 $$E[ f(M_t, R_t)\ ; \  A_t < 2/y]$$
 
 \noindent This establishes the identity (\ref{expect}) of the proposition.
 
\end{proof}

\vspace{3ex}

\begin{remark}  Proposition \ref{integral} is quite similar to Theorem 1 of \cite{WH}.  The authors study Girsanov density processes and develop necessary and sufficient conditions for a Girsanov process to be a martingale.  Theorem 1 shows that, in great generality, the expected value of a Girsanov density equals the tail probability for a certain stopping time.  

Our arguments proving our proposition are similar to those in  \cite{WH}.  In our case, we extend their result to include expected values of  a function of the Brownian motion process multiplied by a Girsanov density.  The choice $f(x,z)\equiv 1$ is a special case of the theorem in  \cite{WH}.  The reader may see this by making the choice
$$X(t) = \frac{ 2M_t}{a+ A_t}$$
as required by Theorem 1.  To define the correct stopping time, one should consider the process

$$\int_0^tY(u)du:=  -2\log(1 - \frac{1}{a} A_t)$$

One may verify directly that $Y(t)$ satisfies the functional equation mentioned in Proposition 1 of \cite{WH}. 
\end{remark}

\vspace{3ex}

\section{The Distribution of $A_t$}

\vspace{3ex}

%the theorem on the distribution begins here

\begin{theorem}\label{aaa} For  $a > 0 $ the distribution of $A_t $ is given by

\begin{equation}\label{dist}
\Pr\{ A_t \leq a\} = e^{-\frac{2}{a}} E[ \exp\big( \frac{2 M_t}{  a+   A_t } \big)]
\end{equation}
\vskip.2in
\noindent Moreover, the random variable $A_t$ has a continuous, positive probability density function $g_t(a)$ which is simply related to the distribution functions of $A_t$ and $\frac{A_t}{M_t}$:

\begin{equation}\label{density}
g_t(a)= \frac{2}{a^2}\Pr\{ A_t \leq a\}  - \frac{2}{a^2}\Pr\{\frac{A_t}{M_t}\leq a\}
\end{equation}

\end{theorem}

\vspace{3ex}

\begin{proof}

The first assertion of the Theorem follows from Proposition \ref{integral} by making the simple choice

$$f(x,z) \equiv 1$$

\noindent Identity (\ref{expect}) becomes in this case

$$\Pr\{  A_t < 2/y \} = E[\exp(\frac{ M_t}{  y^{-1}+   \frac{1}{2} A_t}-y)]$$

$$= E[\exp(\frac{ M_t}{  y^{-1}+   \frac{1}{2} A_t})]e^{-y}$$

\noindent This identity has the rather surprising corollary that  \emph{the expected value above is finite.}  The integrand involves a \emph {double exponential of Browian motion}.  
Using a result of Dufresne,  Corollary 2.3 in  \cite{D}, one can easily show that for  $\beta > \half$

$$E[\exp(\frac{M_t}{y^{-1} + \beta A_t} ) ] < \infty$$

\noindent but here we have shown finiteness for $\beta = \half$.

Notice that the integrand is a monotone function of $y$.    If $y$ varies over some positive interval 
$$(y_0, y_1)$$
then each integrand is dominated by the integrable random variable
$$\exp(\frac{ M_t}{  y_1^{-1}+   \frac{1}{2} A_t})$$
If a sequence $\{y_k\}$ converges to  $\tilde y $ in this interval, we apply the dominated convergence theorem to prove that the expected value converges to its value for $\tilde y $.  Hence, the expected value is a continuous function of $y$.   It follows that the distribution function of $A_t$ is continuous and we may write the identity as

\begin{equation}\label{four}
\Pr\{  A_t \leq 2/y \} = E[\exp(\frac{ M_t}{  y^{-1}+   \frac{1}{2} A_t})]e^{-y}
\end{equation}

\noindent This establishes equation (\ref{dist}).

We show the  probability density  exists by proving that the right hand expression in equation (\ref{dist}) is differentiable w.r.t. $a$. To see this we consider the case of Proposition \ref{integral} for 

$$f(x,z) =x$$

\noindent Identity (\ref{expect}) becomes

$$E[ M_t ; A_t \leq 2/y ] = E[\frac{ M_t}{  (1+   \frac{y}{2} A_t)^2}\exp(\frac{ M_t}{  y^{-1}+   \frac{1}{2} A_t})]e^{-y}$$
$$=y^{-2}E[\frac{ M_t}{  (y^{-1}+   \frac{1}{2} A_t)^2}\exp(\frac{ M_t}{  y^{-1}+   \frac{1}{2} A_t})]
e^{-y}$$

\noindent As in the previous case, the expected value

$$E[\frac{ M_t}{  (y^{-1}+   \frac{1}{2} A_t)^2}\exp(\frac{ M_t}{  y^{-1}+   \frac{1}{2} A_t})]$$

\noindent is necessarily finite and again the integrand is monotone in $y$.  Therefore, the expression is a continuous function of $y$ as we argued in proving (\ref{dist}).  However, this same expected value arises by formally differentiating the expected value

\begin{equation}\label{diff}
 E[\exp(\frac{ M_t}{  y^{-1}+   \frac{1}{2} A_t})]
 \end{equation}
on the r.h. side of  (\ref{dist}) w.r.t. $y$.  Now, as $y$ varies over some positive interval $(y_0, y_1)$ each integrand
$$y^{-2}\frac{ M_t}{  (y^{-1}+   \frac{1}{2} A_t)^2}\exp(\frac{ M_t}{  y^{-1}+   \frac{1}{2} A_t})$$
is dominated by the integrable random variable

$$y_0^{-2}\frac{ M_t}{  (y_1^{-1}+   \frac{1}{2} A_t)^2}\exp(\frac{ M_t}{  y_1^{-1}+   \frac{1}{2} A_t})$$

\noindent So, the $y-$integral from $y_0$ to $y_1$, which equals
$$\exp(\frac{ M_t}{  y_1^{-1}+   \frac{1}{2} A_t}) - \exp(\frac{ M_t}{  y_0^{-1}+   \frac{1}{2} A_t})$$

\noindent is dominated by the product of $y_1 - y_0$ and an integrable function.  Therefore, the difference quotient for the expected value in (\ref{diff}) converges as $y_0 \to y_1 $ and the limit is the entire expression for
\begin{equation}\label{mmm}
-E[ M_t ; A_t \leq 2/y ] e^y
\end{equation}
The same argument applies to the case $y_1 \to y_0$ so that expression (\ref{diff}) has a derivative which equals the expression (\ref{mmm}).
 Since the distribution function is the product of $e^{-y} $ and expression  (\ref{diff}) we conclude that  the distribution function of $A_t$ has a continuous probability density. We differentiate terms in the identity (\ref{four}) to obtain
 
\begin{equation}\label{first}
 -2y^{-2}g_t(2/y) = -\Pr \{A_t \leq 2/y\} +E[ M_t ; A_t \leq 2/y ] 
   \end{equation}
Now, the expected value

$$E[ M_t ; A_t \leq 2/y ] $$
is another distribution function.  Using the change of measure induced by the factor $M_t$, we see that 
$$B_s - s = W_t$$
is a standard Brownian motion and so the expected value equals

$$\Pr\{ \int_0^t \exp (W_s+s-s/2)ds \leq 2/y\}$$
We substitute this expression into equation (\ref{first}) to obtain
$$g_t(a) = \frac{2}{a^2}\big(\Pr \{A_t \leq a\}  -\Pr\{ \int_0^t \exp ( B_s+s/2)ds \leq a\}\big)$$
We see that the expression for $g_t(a)$ is strictly positive since
$$\int_0^t \exp (B_s+s/2)ds $$
is strictly larger than the random variable $A_t$ for each sample path.  Finally we note that the random variable $\frac{A_t}{M_t}$ has the same distribution as the time integral above.

$$\frac{A_t}{M_t} = \int_0^t \exp (B_s - B_t +t/2 -s/2)ds $$
has the same distribution as

$$ \int_0^t \exp (W_{t-s}+t/2 -s/2)ds $$
where $W_s$ denotes a standard Brownian motion; we may change variables in the time integral to obtain the time integral of geometric Brownian motion with positive drift.  This establishes assertion (\ref{density}) of the theorem.

\end{proof}

\begin{remark}  We can rewrite the density formula (\ref{density}) by combining the two probabilities.  The density equals

\begin{equation}\label{pretty}
g_t(a) = \frac{2}{a^2} \Pr \{  \int_0^t \exp (B_s-s/2)ds  \leq a <  \int_0^t \exp (B_s+s/2)ds \}
\end{equation}
\vspace{2ex}

\noindent This expresses the probability density for $A_t$ in terms of a single condition on Brownian motion sample paths up to time $t$.

The existence of a continuous probability density for $A_t$ is a corollary of Proposition 2 of Yor \cite{Y92}.  In the proposition  a conditional density for $A_t$ is given as an integral transform of various transcendental functions.  No explicit connection is made between the density and the distributions of $A_t$ and $A_t / M_t$.  

In Dufresne \cite{D} nice formulas are obtained for the density of a reciprocal of a time integral for some values of a drift parameter in the Brownian motion.  Our choice corresponds to the choice of $\mu = -1$, and density is given as an integral transform in \cite{D}.
  
\end{remark}

\vspace{3ex}

\begin{remark}  The density $g_t(a)$ is a solution of the PDE derived in [BTW, equation (26)]  where it is shown that time integrals of more general geometric Brownian motions have smooth densities.  In particular, the random variable 
$$\frac{A_t}{M_t}$$
has a smooth density since it has the same distribution as the time integral of $\exp (B_s + s/2)$.  The PDE, equation (26), has a simple form in the case of the density of $A_t$:
$$\frac{\partial g}{\partial t} =\frac{\partial ^2 }{\partial a ^2}\big \{\frac{a^2 g}{2} \big \}  - \frac{\partial g}{\partial a}  $$
But, equation (\ref{density}) shows that $a^2 g / 2$ is the difference of two distribution functions.  Therefore, the PDE becomes
$$\frac{\partial g}{\partial t} =\frac{\partial }{\partial a }\big \{g - 
 \frac{\partial }{\partial a }\Pr\{\frac{A_t}{M_t}\leq a\} \big \}  - \frac{\partial g}{\partial a}  $$
$$= - \frac{\partial ^2 }{\partial a ^2}\Pr\{\frac{A_t}{M_t}\leq a\}$$
That is, the time derivative of the density is obtained by differentiating the density of $A_t / M_t$.

\end{remark}
\vspace{3ex}

\begin{corollary}\label{super}  For each $y>0$ the process

$$Z_t = \exp\big( \frac{ M_t}{  y^{-1}+   \frac{1}{2} A_t } \big)$$
is a supermartingale.  In particular,

\begin{equation}\label{yyy}
E[ Z_t] = e^y\Pr\{A_t \leq 2/y\}
\end{equation}
\vskip.1in

\noindent so that $E[ Z_t]$ is a strictly decreasing function of $t$.
\end{corollary}

\vspace{3ex}

\begin{proof}
The process
$$Y_t =  \frac{ M_t}{  y^{-1}+   \frac{1}{2} A_t }$$
also satisfies the SDE that appears in  Lemma \ref{sde}.  A simple calculation shows that 

\begin{equation}\label{dydy}
dY_t =  Y_t dB_t  -  \frac{1}{2} Y_t^2dt
\end{equation}
Therefore, the remark following Lemma \ref{sde} applies to the process $Y_t$: 
\vskip.1in

 Since
$$
Y_t -y=  \int_0^t Y_s dB - \frac{1}{2}  \int_0^tY_s^2ds,
$$
$Y_t-y$ is the exponent of the Girsanov density process which $Y_t$ generates.  A simple stopping time argument, similar to one in the proof of  Proposition \ref{tilde}, shows that
$$Z_t = \exp ( Y_t)$$
is a positive local martingale.  However, the integral identity (\ref{dist}) of the Theorem implies that
$$ E[\exp ( Y_t)]$$
\vskip.1in
\noindent is decaying function of $t$.   Consequently,   $Z_t$ is a local martingale but it is not a martingale.
\end{proof}

\begin{remark}  The process $Y_t$ appears implicitly in Lemma 2.1 of  \cite{BTW}.  In order to derive PDE's for certain expected values involving $A_t$,  the authors compute a diffusion equation for processes slightly more general than

$$ (Y_t^{-1} , M_t).$$
The PDE identities do not apply to equation (\ref{yyy}) since $Z_t$ is not a homogeneous function of $Y_t^{-1}$.

One can derive the corollary from Theorem 1 of \cite{WH}, but the argument here connects the result directly to the behavior of time integrals of geometric Brownian motion.

\end{remark}
\vspace{3ex}

\begin{corollary}\label{moment} 

\begin{equation}\label{expA}
E[\exp\big( \frac{ 2}{A_t} \big)]=\infty
\end{equation}

 \end{corollary}

\vspace{1ex}

\begin{proof} For each $a>0$

$$E[\exp\big( \frac{ 2}{A_t} \big)\ ; A_t \leq a] \geq \exp (2/a)\Pr\{ A_t \leq a\}$$
and equation (\ref{dist}) of the theorem implies that the r.h. expression equals

$$E[ \exp\big( \frac{2 M_t}{  a+   A_t } \big)]  $$
But, this quantity increases as $a\to 0$ and therefore the random variable

$$\exp\big( \frac{ 2}{A_t} \big)$$
is not integrable.
\end{proof} 

\vspace{3ex}

%brownian motions with non-zero drift

\section{Geometric Brownian Motion with Drift}

\vspace{3ex}

\begin{definition} For each $\nu \in \mathbb{R}$  we let $A_t^{(\nu)}$ denote the time integral

\begin{equation}\label{new}
A_t^{(\nu)} = \int_0^t \exp (B_s+\nu s-s/2)ds 
\end{equation}
\end{definition}

\vspace{2ex}
\noindent The random variable $A_t$ in the previous sections is $A_t^{(0)}$ with this notation.

\vspace{3ex}

%the theorem with the drift begins here

\begin{theorem}\label{withdrift} For  $a > 0 $ the distribution of $A_t^{(\nu)}$ is given by

\begin{equation}\label{distdrift}
\Pr\{ A_t^{(\nu)}\leq a\} = a^{2\nu}e^{-\frac{2}{a}} 
E[(a+A_t^{(\nu)})^{-2\nu} \exp\big( \frac{2 \exp (B_t +\nu t-t/2)}{  a+   A_t^{(\nu)} } \big)]
\end{equation}

\end{theorem}

\vspace{3ex}

\begin{proof} For $y >0$ we apply Proposition \ref{integral} to evaluate 

$$E[ ( M_t )^\nu \ ;\ A_t \leq 2/y ]$$
The choice of $f(x,z) = x^\nu$ in the proposition gives the expected value

$$E[ \frac{M_t^\nu}{(1+\frac{y}{2}A_t)^{2\nu}}\exp(\frac{ M_t}{  y^{-1}+   \frac{1}{2} A_t}-y)  ]$$
Next, we multiply these expected values by $\exp(\nu t/2 - \nu^2 t/2)$ so that

$$ M_t^\nu \exp(\nu t/2 - \nu^2 t/2)=  \exp (\nu B_t - \nu^2 t/2)$$
We use this exponential martingale factor to change measure in each integral so that the process

$$\tilde B_s = B_s - \nu s$$
is a standard Brownian motion for $s \leq t$.  We see that 

$$E[ M_t^\nu\exp(\nu t/2 - \nu^2 t/2) \ ;\ A_t \leq 2/y ]
=\Pr\{ A_t^{(\nu)} \leq 2/y\}$$

\noindent while the other expected value equals

$$E[(1+\frac{y}{2}A_t^{(\nu)})^{-2\nu}  \exp(\frac{ \exp (B_t +\nu t - t/2)}{  y^{-1}+   \frac{1}{2} A_t^{(\nu)}}-y)  ]$$
This establishes identity (\ref{distdrift}).

\vspace{3ex}

\end{proof}

\vspace{3ex}

\begin{corollary}\label{tail}  As $a \downarrow 0$ the function

$$\frac{\exp (2/a)}{a}\Pr\{ A_t^{(1/2)}\leq a\}$$
increases.

\end{corollary}

\vspace{3ex}

\begin{proof} For the case that  $\nu = 1/2$,  the formula in (\ref{distdrift}) for the distribution function becomes

\begin{equation}\label{nodrift}
\Pr\{ A_t^{(1/2)}\leq a\} = ae^{-\frac{2}{a}} 
E[\frac{1}{a+A_t^{(1/2)}} \exp\big( \frac{2 \exp (B_t )}{  a+   A_t^{(1/2)} } \big)]
\end{equation}
The integrand of the expected value in (\ref{nodrift}) increases as  $a \downarrow 0$.

\end{proof}
\vspace{3ex}

%the exponential moment of the driftless case is here

\begin{corollary}  The following expected value is infinite.

\begin{equation}\label{classic}
E[\exp\big( \frac{2 }{ \int_0^t \exp (B_s )ds  }\big)] = \infty
\end{equation}
\end{corollary}

\vspace{3ex}

\begin{proof} The corollary states that $E[\exp\big( \frac{2 }{ A_t^{(1/2)} }\big)] = \infty$.  Let $F(a)$ denote the distribution function of $A_t^{(1/2)}$ and consider

$$E[\exp\big( \frac{2 }{ A_t^{(1/2)} } \big)\ ; \ A_t^{(1/2)} \leq 1]$$
We use a standard argument that justifies integration by parts:

$$e^{\frac{2}{a}}= \int_{a}^\infty \frac{2}{x^2}e^{\frac{2}{x}}dx$$
so that

$$E[\exp\big( \frac{2 }{ A_t^{(1/2)} } \big)\ ; \ A_t^{(1/2)} \leq 1]
= \int_{0}^1 \int_{a}^\infty \frac{2}{x^2}e^{\frac{2}{x}}dxdF(a) $$

$$=\int_{0}^\infty \frac{2}{x^2}e^{\frac{2}{x}}  \int_{0}^{x\wedge 1}  dF(a)dx $$
$$=\int_{0}^\infty \frac{2}{x^2}e^{\frac{2}{x}}  F(x\wedge 1)dx $$
But, Corollary \ref{tail} implies that

$$\frac{1}{x^2}e^{\frac{2}{x}}  F(x) \geq \frac{k}{x}$$
 for $x < 1$  where the constant $k > 0$.  Hence, the integral is infinite.

\end{proof}

\vspace{3ex}

%this section has finite exponential moments

\section{Finite Exponential Moments}

\vspace{3ex}

%begin the lemma on A_t (.5)

\begin{lemma}\label{ordinary}  For any $t > 0$
\begin{equation}\label{basic}
E[\exp\big( \frac{1 }{2 \int_0^t \exp (B_s )ds  }\big)] < \infty
\end{equation}

\end{lemma}

\vspace{2ex}

\begin{proof}  The expected value in (\ref{basic}) is

$$E[ \exp\big( \frac{1}{2A_t^{(1/2)} } \big) ]  $$
with the notation of Section 5.  And, an upper bound for the expected value is

$$e+ \int_e^\infty \Pr\{ \exp\big( \frac{1}{2A_t^{(1/2)} } \big) \geq x\} dx$$
Let $x=e^s$ to obtain the following form of the integral above:
\begin{equation}\label{form}
\int_1^\infty \Pr\{  \frac{1}{2A_t^{(1/2)} }  \geq s\}e^s ds
\end{equation}
$$=\int_1^\infty \Pr\{ \  \frac{1}{2s}   \geq A_t^{(1/2)}\ \}e^s ds$$
\noindent It suffices to prove that the integral from $k$ to $\infty$ is finite where $k$ is chosen so that on the interval of integration
$$\frac{4}{s} \leq t$$
For any $s$ in the interval we have

$$A_t^{(1/2)} =   \int_0^t \exp (B_u )du \geq \int_0^{4/s} \exp (B_u )du$$
$$=\frac{4}{s}\cdot \frac{s}{4}\int_0^{4/s} \exp (B_u)du \geq    \frac{4}{s}\exp (\frac{s}{4}
\int_0^{4/s}B_u du )$$
by Jensen's inequality.  Since the random variable

$$ \frac{s}{4}\int_0^{4/s}B_u du $$

\noindent is normal with standard deviation
$$ \frac{2}{\sqrt{3s}} $$
we may replace it with 
$$ \frac{2}{\sqrt{3s}} Z$$
where $Z$ denotes a standard normal random variable.  The integrand in equation (\ref{form}) is dominated by
$$ \Pr\{  \frac{1}{2s}   \geq \frac{4}{s}\exp \left (\frac{2}{\sqrt{3s}} Z\right) \}e^s$$
$$ = \Pr\{  \frac{1}{8}   \geq \exp  \left ( \frac{2}{\sqrt{3s}} Z\right) \}e^s$$
$$ < \Pr\{  -2   \geq \frac{2}{\sqrt{3s}} Z \}e^s$$
since $\log(1/8) < -2$.  We may write this as
$$ \Pr\{  -\sqrt{3s}   >  Z \}e^s$$
Hence the integral in equation (\ref{form}),
$$\int_1^\infty \Pr\{ \  \frac{1}{2s}   \geq A_t^{(1/2)}\ \}e^s ds,$$
 is dominated by

$$e^k+\int_k^\infty \Pr\{  -\sqrt{3s}   >  Z \}e^s ds$$
$$\approx e^k+\int_k^\infty \exp(-\frac{3s}{2})e^s ds < \infty$$

\end{proof}

\vspace{3ex}

%begin the theorem on   M / A

\begin{theorem}\label{MoverA}  For any $\theta < 2$ and $t_0 > 0$ the process

$$\exp\big( \frac{\theta M_t }{A_t} \big) $$

\noindent is a supermartingale for $t \geq t_0$.  In particular,

\begin{equation}\label{m/a}
E[\exp\big( \frac{\theta M_t }{A_t}\big)] < \infty
\end{equation}

\vspace{1ex}
\noindent and  the expected value is a strictly decreasing function of $t$.

\end{theorem}

\vspace{2ex}

\begin{proof} We have seen (proof of Theorem \ref{aaa}) that $A_t / M_t$ has the same distribution as 

$$ \int_0^t \exp (B_s+s/2 )ds  $$

Since this random variable is larger than $A_t^{(1/2)}$, Lemma \ref{ordinary} implies that

$$E[\exp\big( \frac{ M_t }{2A_t}\big)] < \infty$$
Now let 
\begin{equation}\label{uuu}
U_t = \frac{e^{ct}}{4}\frac{M_t}{A_t}
\end{equation}
for any fixed $c > 0$. If $t_0 > 0$ is sufficiently small, so that $ e^{c t_0}/4 \leq 1/2$, we have

$$E[ \exp(U_{t_0})] < \infty$$ 
We define $t_1$ by

$$ e^{c t_1}/4 = \theta$$
and we claim that
$$E[ \exp(U_{t_1})] < \infty$$ 
The integrabliity of $\exp\big( \frac{\theta M_t }{A_t}\big)$ will follow because the choice of $c$ is arbitrary. We first consider the SDE for $U_t$:

$$dU_t = cU_t dt + U_t dB_t - \frac{e^{ct}}{4}\frac{M_t}{A_t^2}M_t dt$$
$$=U_t dB_t + cU_t dt- \frac{4}{e^{ct}}U_t^2dt$$
$$=U_t dB_t + U_t \{c- 4e^{-ct}U_t\}dt$$
Next, we compute the SDE for the process $\exp ( U_t)$.

$$d\exp ( U_t) = U_t \exp ( U_t)dB_t +\exp ( U_t)U_t \{c- 4e^{-ct}U_t\}dt + \frac{1}{2}\exp (U_t)U_t^2dt$$
In order to compute an expected value, we introduce the stopping times

$$\tau_n = \inf \{ t \geq t_0 :  U_t \geq n\}$$
and we obtain
$$E[ \exp ( U_{\tau_n\wedge t_1}) ] - E[ \exp ( U_{ t_0} ) ] $$
$$ = E[ \int_{ t_0}^{\tau_n\wedge t_1}\exp ( U_t)U_t \{c+
 \frac{1}{2}U_t- 4e^{-ct}U_t\} \,dt] $$
$$ \leq E[\int_{ t_0}^{\tau_n\wedge t_1}\exp ( U_t)U_t \{c+
 \frac{1}{2}U_t- 4e^{-ct_1}U_t\} \,dt] $$
\begin{equation}\label{difference}
 = E[ \int_{ t_0}^{\tau_n\wedge t_1}\exp ( U_t)U_t \{c+
 \frac{1}{2}U_t- \theta^{-1}U_t\} \,dt] 
 \end{equation}
Notice that  if
$$U_t \geq b$$
where 
$$b :=c(\theta^{-1} - \frac{1}{2} )^{-1} $$
then the integrand in equation (\ref{difference}) is negative.
It follows that 
$$E[ \exp ( U_{\tau_n\wedge t_1}) ] - E[ \exp ( U_{ t_0}) ] $$
$$ \leq E[ \int_{t_0}^{\tau_n\wedge t_1}\exp ( U_t)U_t \{c+
 \frac{1}{2}U_t- \theta^{-1}U_t\}1_{\{ U_t \leq b\}}dt] $$
 Now we take the limit as $n\to \infty$  and, using the condition that $U_t \leq b$ in the integrand,  we see that

  $$E[ \exp ( U_{ t_1} )]  < \infty$$
\vskip.1in
This establishes that the expected value in equation (\ref{m/a}) is  finite.  It is a  strictly decreasing function of $t$ because the random variable $\frac{ A_t }{M_t}$ has the same distribution as $A_t^{(1)}$ which is a strictly increasing function of $t$.  It remains to establish the supermartingale property.  Corollary \ref{super} implies that for any $\alpha < 1$ the process 

$$\big\{\exp\big( \frac{2 M_t}{  a+   A_t } \big)\big\}^\alpha$$
is a supermartingale because it is a concave function of a supermartingale.  Now as $a\to 0 $ the pointwise limit of this process is 

$$\big\{\exp\big( \frac{2 M_t}{  A_t } \big)\big\}^\alpha
 = \exp\big( \frac{2\alpha M_t}{  A_t } \big)$$
That is, the process in the statement of the theorem is integrable and is the limit of non-negative supermartingales.  Hence, it is also a supermartingale.
\end{proof}

\vspace{3ex}

\begin{remark}  It follows from Theorem \ref{withdrift} that the distribution function of $\frac{ M_t }{A_t}$
has the form
$$a^2 \exp ( -2/a)K_a$$
\vspace{1ex}

\noindent where the function $K_a$ increases as $a\downarrow 0$.  So, if $\theta > 2$ in equation (\ref{m/a}), the expected value is infinite.  It is unclear for the case $\theta = 2$ if the expected value is finite.

\end{remark}

\vspace{3ex}

%begin the corollary about A_t (.5)

\begin{corollary}   For any $\theta < 2$ and $t > 0$

\begin{equation}\label{finite}
E[\exp\big( \frac{\theta }{ \int_0^t \exp (B_s )ds  }\big)] < \infty
\end{equation}

\end{corollary}

\vspace{2ex}

\begin{proof}  Since the expected value in (\ref{finite}) is a decreasing function of $t$, it suffices to show the expected value is finite for $t$ arbitrarily small.  Since $\frac{M_t}{A_t}$ has the same distribution as

$$ \int_0^t \exp (B_s+s/2 )ds,  $$

\noindent Theorem \ref{MoverA} implies that for any $\tilde \theta < 2$
$$E[\exp\big( \frac{\tilde \theta  }{\int_0^t \exp (B_s+s/2 )ds}\big)] < \infty$$
For a given $\theta < 2$, choose $t$ so small that
$$\tilde \theta :=\theta \exp(\frac{t}{2}) < 2$$
so that 
$$E[\exp\big( \frac{ \theta  }{e^{-t/2}\int_0^t \exp (B_s+s/2 )ds}\big)] < \infty$$
We see that the expected value above is larger than
$$E[\exp\big( \frac{\theta }{ \int_0^t \exp (B_s )ds  }\big)] $$
and so this expected value is finite.\end{proof}

\vspace{3ex}

\begin{remark}  Some related work on exponential moments of $A_t$ is mentioned in Yor \cite{Y92}.
Equation (1.e) of the article states that
$$E[\ \frac{1}{\sqrt{A_t}}\exp\big(-\frac{u^2}{2A_t}\big)\ ] = \frac{1}{\sqrt{(1+u^2)t}}\exp\big(
 -\frac{1}{2t}(\sinh^{-1} u)^2\big)$$
However, a difference in notation requires that $A_t$ be expressed with our notation as 
$\frac{1}{4}A_{t/4}^{(1/2)}$ .  By writing $t' = t/4$ we see that the formula in Equation (1.e) is an expression for
$$E[\ \frac{2}{\sqrt{A_{t'}^{(1/2)}}}\exp\big(-\frac{2u^2}{A_{t'}^{(1/2)}}\big)\ ] $$
Equation (\ref{finite})  implies that the expected value is finite for all complex values of $u$ such that $Re(u^2) >-1$ and is analytic on this region.  Since
$$\sinh^{-1} u = \log(u+\sqrt{1+u^2})$$
one sees that the right hand expression also has an analytic extension.  The formula has a singularity at the value $u=i$ which corresponds to the infinite exponential moment of Corollary \ref{classic}.

In addition, Theorem 4.1 of \cite{D} provides an integral formula for the density of $A_t^{(1/2)}$ and one can show that the expected value is finite (for $\theta < 2$) using Theorem 4.1.  The corollary provides a probabilistic  proof of this fact.  The result that the moment is infinite for $\theta =2$ is new.
\end{remark}

\vspace{3ex}

An expected value considered in \cite{Y92}, \cite{RS}, and in other works concerning time integrals in financial mathematics is the `price' for an Asian call option.  We present a new expected value for the simplest option price and indicate how to derive a corresponding  formula for arbitrary constant drift and volatility values.

\vspace{3ex}

\begin{proposition}  For any $a>0$

\begin{equation}\label{call}
E[\ (A_t - a)^{+}\ ] = t-a +a^2 E[(a+A_t)^{-1} \exp\big( \frac{2 M_t}{  a+   A_t } -\frac{2}{a} \big)]
\end{equation}
\end{proposition}

\vspace{3ex}

\begin{proof}  The notation $(A_t - a)^{+}$ involves the indicator function of the event $A_t \geq a$ so it is natural to consider
$$ E[\ A_t - 2/y \ ; \ A_t < 2/y \ ] $$
$$=-\frac{2}{y}E[\  1 - \frac{y}{2}A_t  \ ; \ A_t < 2/y \ ] $$
We apply Proposition \ref{integral} where we make the choice
$$f(x,z) = \frac{x}{-z}$$

Then
$$-2f(M_t, R_t) = -2M_t \frac{y^{-1} - \frac{1}{2}A_t}{M_t}$$
$$= -\frac{2}{y}(1 - \frac{y}{2}A_t)$$
Equation (\ref{expect}) shows that the expected value over $A_t < 2/y$ is given by

$$-\frac{2}{y}E[ (1+ \frac{y}{2}A_t)^{-1}\exp(\frac{ M_t}{  y^{-1}+   \frac{1}{2} A_t}-y)] $$
On the other hand, 
$$E[\ A_t - \frac{2}{y}\ ] = t - \frac{2}{y}$$
so we take the difference of these expected values and let $\frac{2}{y} = a$ to obtain the result.
\end{proof}

\vspace{3ex}

\begin{remark}
The proposition requires the Brownian motion to have  drift $-1/2$.  To derive an expectation formula for arbitrary drift  one can apply Proposition \ref{integral} to the expected value
$$ E[(M_t)^\nu (\ A_t - 2/y )\ ; \ A_t < 2/y \ ] $$
The factor of $(M_t)^\nu$ is relevant for a change of measure (see the proof of Theorem \ref{withdrift}) so that the time integral will contain an arbitrary drift.  To incorporate a volatility factor in the Brownian motion, one can make a simple time scale change.  We let $s=\sigma^2 s'$ so that
$$\frac{1}{\sigma^2}\int_0^t \exp (B_s)ds = \int_0^{t'} \exp (B_{\sigma^2 s'})ds'$$
This random variable has the same distribution as the time integral of $\exp (\sigma B_{ s})$.
\end{remark}


\vspace{3ex}



\begin{thebibliography}{DGMS}

\bibitem[BY]{BY} J. Bertoin and M. Yor: \emph{Exponential functionals of L\' evy processes}, Probabilty Surveys \textbf{2}, 191-212 (2005).

\bibitem[BTW]{BTW} R. Bhattacharya, E.Thomann, and E. Waymire: \emph{A note on the distribution of integrals of geometric Brownian motion}, Stat. and Prob. Letters  \textbf{55}, 187-192 (2001).

\bibitem[D]{D} D. Dufresne: \emph{The integral of geometric Brownian motion}, Adv. in Appl. Probab. \textbf{33}, 223-241 (2001).


\bibitem[D2]{D2} D. Dufresne: \emph{The distribution of a perpetuity, with applications to risk theory and pension funding}, Scand. Actuar. I. \textbf{1-2}, 39-79 (1990).


\bibitem[GY]{GY}H. Geman and M. Yor: \emph{Asian Options, Bessel Processes and Perpetuities}, Math. Finance \textbf{2}, 349-375 (1993).

\bibitem[KS]{KS} I. Karatzas and S.Shreve: \emph{Brownian Motion
and Stochastic Calculus}, Springer-Verlag New York (1991).

\bibitem[K]{K} K. Kim: \emph{Moment Generating function of the inverse of integral of geometric Brownian Motion}, Proc. Amer. Math. Soc. \textbf{132}, 2753-2759 (2004)

\bibitem[MY1]{MY1}H. Matsumoto and M. Yor: \emph{An analogue of Pitman's 2M -- X theorem for exponential Wiener functionals, Part I:  A time-inversion approach.}, Nagoya Math. J. \textbf{159}, 125-166  (2000).

\bibitem[MY2]{MY2}H. Matsumoto and M. Yor: \emph{An analogue of Pitman's 2M --- X theorem for exponential Wiener functionals, Part II:  The role of the generalized inverse Gaussian laws.}, Nagoya Math. J. \textbf{162}, 65-86  (2001).

\bibitem[MY3]{MY3}H. Matsumoto and M. Yor: \emph{A relationship between Brownian motions with opposite drifts via certain enlargments of the Brownian filtration}, Osaka Journal of Mathematics \textbf{38}, 383-398 (2001).

\bibitem[RS]{RS} L.C.G. Rogers and Z. Shi: \emph{The value of an Asian option}, J. Appl Appl. Probab.\textbf{32}, 1077-1088 (1995).

\bibitem[WH]{WH} B. Wong and C.C. Heyde: \emph{On the martingale property of stochastic exponentials}, J. Appl. Probab. \textbf{41}, 654-664 (2004).


\bibitem[Y(92)]{Y92} M. Yor: \emph{On some exponential functionals of Brownian
motion}, Adv. in Appl. Probab. \textbf{24}, 509-531 (1992).

\bibitem[Y01]{Y01} M. Yor: \emph{Exponential Functionals of Brownian Motion and Related Processes.}, Springer-Verlag, New York. (2001).

\end{thebibliography}
\end{document}